\documentclass[12pt,leqno]{article}
\usepackage{amssymb, amscd, amsmath, amsthm}

\allowdisplaybreaks
\def\thefootnote{}

\def\beq{\begin{equation}}
\def\eeq{\end{equation}}

\theoremstyle{plain}
\newtheorem{Thm}{Theorem}

\theoremstyle{definition}
\newtheorem*{rema}{Remark}

\numberwithin{equation}{section}
\title{The exceptional set of Goldbach numbers within the set of squares}

\author{Alberto Perelli and J\'anos Pintz\thanks{Supported by the National Research Development and Innovation Office, NKFIH,  K~133819, K~147153 and Excellence 151341.}}


\footnotetext{ Exceptional Goldbach numbers among squares.}

\begin{document}
\maketitle
\renewcommand{\thefootnote}{\arabic{footnote}}
\setcounter{footnote}{0}

\date{}

\begin{abstract}
We prove that the set of numbers of type $n = 2m^2$ below $X$ which cannot be written as a sum of two primes has cardinality at most $X^{19/24 + o(1)}$.
\end{abstract}

\maketitle
\section{Introduction}
\label{sec:1}

The second named author proved very recently in a joint work with I. Alth\"ofer \cite{PA2026} that the exceptional set in Goldbach's problem, i.e.\ ($\mathcal P$ denotes the set of primes)
\beq
\label{eq:1.1}
E(X): \# \bigl\{n; 2\mid n, n\leq X \ \ n \neq p_1 + p_2 \ \text{ if } p_i \in \mathcal P\bigr\}
\eeq
satisfies unconditionally
\beq
\label{eq:1.2}
E(X) \ll_\varepsilon X^{7/12 + \varepsilon}\ \text{ for any } \varepsilon > 0.
\eeq
We note that under the Generalized Riemann Hypothesis (GRH) for Dirichlet's $L$-functions, the classical result of Hardy and Littlewood \cite{HL1924} was
\beq
\label{eq:1.3}
E(X) \ll_\varepsilon X^{1/2 + \varepsilon}\ \text{ for any } \varepsilon > 0.
\eeq
Clearly, even \eqref{eq:1.3} gives no nontrivial result for
\beq
\label{eq:1.4}
E_{k}(N) = \# \left\{n, n \leq N, \ 2n^k \neq p_1 + p_2\ \text{ if } p_i \in \mathcal P\right\}
\eeq
for any $k \geq 2$.

\smallskip
The first named author proved in 1996 \cite{Per1996} for any $k$ and $A$
\beq
\label{eq:1.5}
E_k(N) \ll_{k, A} N(\log N)^{-A}.
\eeq
This was improved by Br\"udern, Kawada and Wooley \cite{BKW2001} to
\beq
\label{eq:1.6}
E_k(N) \ll_k N^{1 - c_0/k}
\eeq
with an unspecified absolute constant $c_0 > 0$. The above works handled also the case of a general polynomial $\Phi(x) \in \mathbb Z[x)$ in place of $x^k$.

\smallskip
Our goal is to show a stronger explicit estimate for $k = 2$ (and later for larger values of $k$).

\hfill \newpage
\begin{Thm}
\label{th:1}
$E_2(N) \ll_\varepsilon N^{19/24 + \varepsilon}$.
\end{Thm}

We will use the circle method $(\varepsilon > 0$ arbitrary) with
\beq
\label{eq:1.7}
S(\alpha) := \sum_{X_1 < p \leq X} \log p e (p\alpha), \ \ X_1 = X^{1 - \varepsilon}, \ e(u) = e^{2\pi i u}
\eeq
and
\beq
\label{eq:1.8}
f(\alpha) := \sum_{1 \leq n \leq N} e(2\alpha n^2), \ \ N = \sqrt{X/2}.
\eeq
The Farey arcs are defined with a parameter $P_0 \leq \sqrt{X}$ and with $Q_0 = X/P_0$ which lead to the following major and minor arcs:
\beq
\label{eq:1.9}
\mathfrak M = \bigcup_{q \leq P_0} \bigcup_{\substack{a\\ (a, q) = 1}} \left[\frac a{q} - \frac 1{qQ_0}, \frac{a}{q} + \frac1{q Q_0} \right], \ \ \mathfrak m = \left[\frac1{Q_0}, 1 + \frac1{Q_0}\right] \setminus \mathfrak M.
\eeq

Let $m \in [X/2, X]$. Our goal is to reach
\beq
\label{eq:1.10}
\mathcal R(m) = \sum_{\substack{p + p' = m\\ p, p' \in [X_1, X]}} \log p \cdot \log p' > 0
\eeq
for values of type $m = 2n^2$, apart from a small exceptional set.
We consider first this for all $m \in [X/2, X]$.
The contribution of the major and minor arcs, respectively, is
\beq
\label{eq:1.11}
R_1(m) = \int\limits_{\mathfrak M} S^2(\alpha) e(-m\alpha)d\alpha, \ \ \
R_2(m) = \int\limits_{\mathfrak m} S^2(\alpha) e(-m\alpha)d\alpha,
\eeq
while
\beq
\label{eq:1.12}
R(m) = \int\limits_0^1 S^2(\alpha) e(-m\alpha)d\alpha = R_1(m) + R_2(m).
\eeq
In case of the unrestricted binary Goldbach problem our goal is, in general, to show, with a possibly large value of the parameter $P_0$, an asymptotic evaluation of type
\beq
\label{eq:1.13}
R_1(m) = (1 + o(1)) \mathfrak S(m),
\eeq
where
\beq
\label{eq:1.14}
\!\!\!\mathfrak S(m)\! =\! 2C_0 \prod_{\substack{p\mid m\\ p > 2}} \Bigl(\!1\! +\! \frac1{p\! -\! 2}\Bigr), \ \
C_0\!  =\! \prod_{p > 2} \Bigl(\!1\! -\! \frac1{(p\! -\! 1)^2}\Bigr)\! = 0.66016...
\eeq
and to estimate from above the average of $\left|R_2(m)\right|$. In practice, in most cases, instead of \eqref{eq:1.13}, it is sufficient to show
\beq
\label{eq:1.15}
R_1(m) \geq c \mathfrak S(m),
\eeq
where (and in the following) $c$ is a generic absolute constant.

\section{Auxiliary results}
\label{sec:2}

The second named author \cite{Pin2023} succeeded in an evaluation of $R_1(m)$ for a parameter value as large as
\beq
\label{eq:2.1}
P_0 = X^{4/9 - \varepsilon} \ \ \text{ for any } \ \varepsilon > 0.
\eeq
However, a relatively small set $\mathcal E(H, P_0, X)$ of possible ``generalized exceptional zeros'' of primitive Dirichlet $L$-functions (with $\text{\rm cond }\chi \leq P_0$) appeared in the formula.
They were defined as ($L = \log X$, $H = o(L)$)
\beq
\label{eq:2.2}
\aligned
\mathcal E(H, P_0, X) = &\biggl\{\!(\varrho, \chi);\, \varrho = \beta + i\gamma;\, \exists \chi, \, \text{\rm cond}\, \chi \leq P_0, \, L(\varrho, \chi) = 0, \\
&\qquad \beta \geq 1 - \frac{H}{L}, \, |\gamma|\!\leq\! \sqrt{X}\!\biggr\}
\endaligned
\eeq
with additional generalized singular series $\mathfrak S(x_1, x_2, m)$, $\mathfrak S(x, m)$ satisfying
\beq
\label{eq:2.3}
\left|\mathfrak S(\chi_1, \chi_2, m)\right| \leq \mathfrak S(m), \ \ |\mathfrak S(\chi, m)| \leq \mathfrak S(m).
\eeq

\begin{rema}
We can choose $H = o(L)$ with $H \to \infty$ when $X \to \infty$ (e.g.\ $H = \sqrt{L}$).
\end{rema}

\medskip
\noindent
{\bf Theorem A.} {\it For every $P \leq X^{4/9 - \varepsilon}$ we can choose $P_0 \in \bigl[PX^{-\varepsilon}, P\bigr]$ such that with the notation above and $\mathcal E = \mathcal E(H, P_0, X)$
we have}
\beq
\label{eq:2.4}
\aligned
R_1(m) &= \mathfrak S(m) + \sum_{(\varrho_i, \chi_i) \in \mathcal E}
\sum_{(\varrho_j, \chi_j) \in \mathcal E} \mathfrak S(\chi_i, \chi_j, m) \frac{\Gamma(\varrho_i)\Gamma(\varrho_j)}{\Gamma(\varrho_i + \varrho_j)} m^{\varrho_i + \varrho_j - 1}\\
&\qquad -\sum_{(\varrho_i, \chi_i) \in \mathcal E} \mathfrak S(\chi_i, m) \frac{m^{\varrho_i}}{\varrho_i} + O_{\varepsilon} \left(\mathfrak S(m) Xe^{-c_\varepsilon H}\right) + O_\varepsilon \left(X^{1 - \varepsilon}\right).
\endaligned
\eeq

This is Theorem 1 of \cite{Pin2023}.

\smallskip
If the set $\mathcal E = \mathcal E(H, P_0, X)$ is empty, then \eqref{eq:2.4} gives the simple asymptotic $\mathfrak S(m)m$ for the contribution of the major arcs.
The very recently published fantastic result of Open AI \cite{OAI2026a}, \cite{OAI2026b} is 

\medskip
\noindent
{\bf Theorem B.} {\it There are no zeros of any Dirichlet $L$-functions (including the Riemann zeta function) in the halfplane $\text{\rm Re } s > 7/8$.}

\medskip
Together with Theorem A this implies

\medskip
\noindent
{\bf Theorem C.} {\it With the above notations we have}
\beq
\label{eq:2.5}
R_1(m) = (1 + o(1)) \mathfrak S(m).
\eeq

\medskip
On the minor arcs one can use the 90-year-old estimate of Vinogradov \cite{Vin1937}, which was proved later in a much simpler form by Vaughan \cite{Vau1977}.
However, again very recently, Maynard, Pandey and Radziwill \cite{MPR2026} gave a striking improvement as

\medskip
\noindent
{\bf Theorem D.} {\it Let $\alpha = a/q + \delta$ with $(a, q) = 1$, $q \leq \sqrt{N}$, $\delta \leq 1/\bigl(qN^{1/2}\bigr)$.
Then with $B = \max(q, qN\delta)$ we have}
\beq
\label{eq:2.6}
\biggl|\sum_{p \leq N} \log p e(p\alpha)\biggr| \leq N^{o(1)} \left(\frac{N}{B^{1/2}} + N^{19/24}\right).
\eeq

\medskip
Finally we will make use of the classical result originating from Fermat and Euler which we will use in a simpler form as

\medskip
\noindent
{\bf Lemma E.}
{\it $a_n := \# \left\{(u, v),\, n = u^2 + v^2,\, u,v \in \mathbb Z^+\right\} \leq 4d(n)$,
where $d(n)$ is the number of divisors of $n$.}

\medskip
Further, we will use the also classical

\medskip
\noindent
{\bf Lemma F.} $ \sum\limits_{n \leq X} d^2(n) \ll X \log^3 X$.

\section{Proof of Theorem \ref{th:1}}
\label{sec:3}

Choosing the Farey arcs and the parameters $P, P_0$ as in Section \ref{sec:1}, with the corresponding $Q_0 = \dfrac{X}{P_0}$, we have for $N/\sqrt{2} \leq m \leq N = \sqrt{X/2}$ by Theorem~C
\beq
\label{eq:3.1}
R_1(2m^2) = (1 + o(1)) \mathfrak S(2m^2) \cdot 2m^2 > 2m^2 \geq N^2.
\eeq
On the other hand, Lemmas E and F imply by Parseval identity
\beq
\label{eq:3.2}
\int\limits_0^1 |f(\alpha)|^4 d\alpha \ll X \log^3 X \ll N^2 \log^3 N.
\eeq

The conditions of Theorem D are satisfied and $\sqrt{P_0} > X^{5/24}$.
This can be seen in the following way.
Let us use Dirichlet's approximation theorem with the parameter $U = X$ and $U = X_1$ in place of $N$ for any $\alpha$ as follows: Let
\beq
\label{eq:3.3}
\alpha = \frac{a}{q} + \delta, \ \ q\leq \sqrt{U}, \ \ |\delta| \leq \frac1{q\sqrt{U}}, \ \ B(U) = \max(q, qU\delta)
\eeq
for an $\alpha \in \mathfrak m$. If $q > P_0$, then clearly $B(U) \geq q > P_0 \geq P_0 \dfrac{U}{X}$.
If $q \leq P_0$, then $\alpha \in \mathfrak m$ implies $\delta > \dfrac1{qQ_0}$, so $B(U) \geq \dfrac{U}{Q_0} = \dfrac{UP_0}{X}$. Thus in both cases ($U = X$ or $X_1$)
\beq
\label{eq:3.4}
\frac{U}{\sqrt{B(U)}} \leq \frac{\sqrt{U}\sqrt{X}}{\sqrt{P_0}} \leq \frac{X}{\sqrt{P_0}} \leq U^{{19/24} + \varepsilon} \leq X^{19/24 + \varepsilon}.
\eeq
Therefore we have by Theorem D
\beq
\label{eq:3.5}
\sup_{\alpha \in \mathfrak m} |S(\alpha)| \leq X^{\frac{19}{24} + \varepsilon}.
\eeq

In connection with the estimate on the minor arcs we will use a neat idea of \cite{BKW2001}.
Since $R_2(m)$ is real for any $m$ let us define
\beq
\label{eq:3.6}
\eta(\ell) = \begin{cases}
1 &\text{ if } R_2(\ell) \geq 0,\\
-1 &\text{ if } R_2(\ell) < 0.
\end{cases}
\eeq
Then
\beq
\label{eq:3.7}
\sum_{m \leq N} \left|R_2(2m^2)\right| = \int\limits_{\mathfrak m} S^2(\alpha) K(-\alpha),
\eeq
where
\beq
\label{eq:3.8}
K(\alpha) = \sum_{m \leq N} \eta (2m^2)e(2m^2 \alpha) .
\eeq

If $t$ is an even number, then $\int\limits_0^1|K(\alpha)|^t d\alpha$ counts the number of solutions of an underlying diophantine equation, with each solution counted with weight at most $1$ in modulus.
Hence,
\beq
\label{eq:3.9}
\int\limits_0^1 |K(\alpha)|^t d\alpha \leq \int\limits_0^1 \Bigl|\sum_{m \leq N} e(2m^2\alpha)\Bigr|^t d\alpha.
\eeq
Consequently, by \eqref{eq:3.2}, \eqref{eq:3.6}--\eqref{eq:3.9} and H\"older's inequality
\beq
\label{eq:3.10}
\aligned
\sum_{m\leq N} \bigl|R_2(2 m^2)\bigr| &\leq \int\limits_{\mathfrak m} |S^2(\alpha)K(\alpha)|d\alpha \\
&\ll \sup_{\alpha \in \mathfrak m} |S(\alpha)|^{1/2} \biggl(\int\limits_0^1 |S(\alpha)|^2 d\alpha\biggr)^{3/4} \biggl(\int\limits_0^1 |f(\alpha)|^4 d\alpha\biggr)^{1/4}\\
&\ll_\varepsilon X^{19/48 + \varepsilon} \cdot X^{3/4} \cdot X^{1/4} \cdot L^c \ll_\varepsilon X \cdot N^{19/24 + \varepsilon}.
\endaligned
\eeq
This means that by \eqref{eq:3.1}
\beq
\label{eq:3.11}
R_1(2m^2) = (1 + o(1))\mathfrak S(2m^2)2m^2 > \bigl|R_2(2m)\bigr|^2,
\eeq
consequently
\beq
\label{eq:3.12}
R(2m^2) > 0
\eeq
apart from an exceptional set of size $N^{19/24 + o(1)}$.\hfill $\square$

\medskip
\noindent
Alberto Perelli, Dipartimento di Matematica, Universit\`a di Genova, Via Dodecaneso 35, 16146 Genova, Italy, e-mail: alberto.perelli@unige.it

\smallskip
\noindent
J\'anos Pintz, HUN-REN Alfr\'ed R\'enyi Institute of Mathematics, Budapest, Re\'altanoda u. 13--15, H-1053 Hungary,
e-mail: pintz@renyi.hu

\end{document}